\documentclass{ws-ijtaf}
\usepackage{array}
\usepackage{amsmath}
\usepackage{amsfonts}
\usepackage{graphicx}
\usepackage{subfigure}
\begin{document}
\title{An Algorithmic Approach to Non-self-financing Hedging in a Discrete-Time Incomplete Market}
\author{N. Josephy, L. Kimball, V. Steblovskaya}
\address{{Mathematical Sciences, Bentley College, 175 Forest Street\\
Waltham, MA 02452-4705, USA }\\ \email{njosephy@bentley.edu,
lkimball@bentley.edu, vsteblovskay@bentley.edu}}
\author{A. Nagaev}
\address{Nicolaus Copernicus University, Torun, Poland}
\author{M. Pasniewski}
\address{ ABS
Analytics Bank Of America, New York, New York,\\
\email{Michal.Pasniewski@bankofamerica.com}}

\newcommand{\cor}[1]{#1}
\newcommand{\msg}[1]{}
\renewcommand{\thesubfigure}{}

\maketitle

\maketitle
 \begin{abstract}
 We present an algorithm producing a dynamic non-self-financing
 hedging strategy in an incomplete market corresponding to investor-relevant risk
 criterion.
 The optimization is a
two stage process that first determines admissible model
parameters that correspond to the market price of the option being
hedged. The second stage applies various merit functions to
bootstrapped samples of model residuals to choose an optimal set
of model parameters from the admissible set. Results are presented
for options traded on the New York Stock Exchange.
 \end{abstract}

\section{Introduction}
Pricing and hedging of financial assets in incomplete markets is an
active research area in mathematical finance.  One of the possible
ways to produce a model of an incomplete market in discrete time
setting is to assume that stock price relative changes (jumps) can
take more that two values as opposed to the classical binomial model
of Cox-Ross-Rubinstein. In \cite{Tessitore1996}, \cite{wolcyznska},
and \cite{hammarlid} an incomplete market in which stock price jumps
follow a multinomial distribution is studied. In
\cite{Tessitore1996} the no-arbitrage option price interval is
studied, \cite{wolcyznska} and \cite{hammarlid} discuss risk
minimization aspects in option pricing and hedging.

The multinomial model has been further extended to the case where
stock price jumps are distributed over a bounded interval and
options under consideration have convex pay-off functions. In
\cite{Ruschendorf2002}, the upper and lower bounds for no-arbitrage
prices of a European contingent claim with convex pay-off are
obtained. The series of works by A. Nagaev et al. (see
\cite{Nagaev2003}, \cite{WiennaMarkov}, \cite{WiennaNon-Identical},
and \cite{SashaRussian}) are devoted to asymptotic behavior of the
residual value of a minimum cost super-hedge. The residual value
occurs as a result of non-self-financing dynamic hedging strategy of
an option seller introduced and discussed in \cite{Nagaev2003}.

A significant proportion of research on option pricing and hedging
in incomplete markets constructs self-financing trading strategies
that satisfy both a primary no-arbitrage condition and secondary
conditions on portfolio risk and return. A comprehensive survey of
modern methodologies can be found in \cite{Staum}. A number of
articles that deal with frictions in markets, shortfall risks and
quadratic hedging (all producing incomplete markets) can be found in
the recent compendium \cite{Jouini2001}.

 Less prevalent is the study of non-self-financing trading
strategies in similar economic environments. The encyclopedic
reference \cite{Shiryaev1999} and the more modest
\cite{Melnikov2002} both illuminate option pricing with consumption,
the model which is similar to the work presented here. Our work is
an initial investigation in the algorithmic study of
non-self-financing strategies discussed in \cite{Nagaev2003} and
\cite{SashaRussian}. We explore the short term behavior of the
residual value of a minimum cost super-hedge, whose long term
behavior was studied by Nagaev et al.

Assuming independent and identically distributed jumps in the
underlying stock process, we use historical data and a bootstrap
simulation process to develop an algorithm producing a dynamic
non-self financing hedging strategy. The resulting hedging strategy
constructs a residual sequence with improved investor risk criteria
as compared to other possible hedging strategies. No additional
assumptions are placed on the underlying stock price jump process
other than having bounded support. An algorithmic approach similar
to our use of bootstrap simulation, but having a different
theoretical foundation and goals can be found in
\cite{Sarykalim2004}.

The remainder of the paper is organized as follows. We develop the
discrete time financial model in section \ref{sec:discrete}.  The
notion of the residual value of a minimum cost super-hedge is
developed in sections \ref{sec:superhedge} through
\ref{section:choice}. The algorithm is described in section
\ref{sec:riskmin}. Algorithm implementation and illustrative results
are presented in \ref{sec:algandresults}. We conclude with some
remarks in section \ref{sec:conclusion} and directions for further
study in section \ref{sec:future}.
\section{Discrete Time Financial Model}
\label{sec:discrete}
 Following the theoretical development in
\cite{Nagaev2003}, our discrete time financial model consists of two
fundamental assets and a derivative security.
\begin{enumerate}
\item A risk-free bond with fixed interest rate $r$, evolving
from an initial value $b_0 > 0$ at time $t = 0$ to $b_k$ at time $t
= k$ as
\[
b_k = b_0 (1 + r)^k
\]
\item A risky stock evolving from an initial value $s_0$ at
time $t = 0$ to $s_k$ at time $t = k$ as
\[
s_k = s_0 \xi_1 \xi_2\cdots \xi_k
\]
where $\xi_k = \frac{s_k}{s_{k-1}}$ are assumed to be independent
and identically distributed random variables with probability
distribution having support equal to a bounded interval $[D,U]$. No
further assumptions are made on the distribution function for the
$\xi_k$. However, this assumption is sufficient to render our
financial market incomplete.
\item A derivative security with convex payoff function $f$. For our
numerical investigations, we use a European call option on the
stock. We take the position of an option seller who wishes to hedge
the potential liability of the sold derivative being exercised.
\end{enumerate}

Our goal is to develop and evaluate an algorithm that will determine
a dynamic, non-self-financing hedging strategy consisting of a
portfolio of our stock and bond assets in our incomplete market. The
portfolio will approximately hedge the derivative security and will
satisfy additional criteria, based on the deviation of the portfolio
value from the required hedging value, that are meaningful to the
investor.

\section{Super-hedging Portfolio}
\label{sec:superhedge}
 Based on a convexity argument,
\cite{Nagaev2003} showed that when the parameters $D$ and $U$ are
known,  there is a minimum cost super-hedge whose value at every
time instant $t=k$ is greater than or equal to the value of the
derivative security. This super-hedge is a portfolio of $\gamma_k$
stocks and $\beta_k$ bonds at time $t=k$ described by
\begin{equation}
 \gamma_k(U,D)= \frac{g_{k+1}(U,D,s_kU) - g_{k+1}(U,D,s_kD)}{s_k\left(U - D\right)}
    \label{eq_gamma}
\end{equation}
\begin{equation}
    \beta_k(U,D)= \frac{U g_{k+1}(U,D,s_kD) - D g_{k+1}(U,D,s_kU)}{(1+r) b_k (U - D)}
    \label{eq_beta}
\end{equation}
where
\begin{equation}
 g_k(U,D,s) = (1+r)^{-(n-k)}\sum_{j=0}^{n-k}C_{n-k}^j[p(U,D)]^j
[1-p(U,D)]^{n-k-j}f(sU^jD^{n-k-j}) \label{CRR_price}
\end{equation}
\begin{equation}
p(U,D) = \frac{(1+r) - D}{U - D} \label{eq_pval}
\end{equation}
$n$ is the number of periods to expiration,  $C_{n-k}^j$ is the
binomial coefficient, and $f$ is a convex pay-off function. At each
time instant $k$, the dynamic super-hedge portfolio constructed in
the prior period is liquidated and the proceeds are used to
construct a new portfolio for the current period. The liquidation
value of the prior period portfolio is given by
\begin{eqnarray}
v_k(U,D) &=&\gamma_{k-1}(U,D)s_k+\beta_{k-1}(U,D)b_k \nonumber\\
 & =&\frac{U-\xi
_k}{U-D}g_k(U,D,s_{k-1}D)+\frac{\xi _k-D}{U-D}g_k(U,D,s_{k-1}U)
\label{eq_liquid}
\end{eqnarray}
The funds required to construct the new period portfolio, or set-up
cost, is given by
\begin{equation}
\overline{X}_k(U,D)= g_k(U,D,s_{k-1}\xi _k) \label{eq_setup}
\end{equation}
(see section 4 for more detail). The liquidation value
(\ref{eq_liquid}) will exceed the set-up cost (\ref{eq_setup})
producing a residual amount $\delta_k$
\begin{equation}
\delta _k(U,D)=\frac{U-\xi _k}{U-D}g_k(U,D,s_{k-1}D)+\frac{\xi
_k-D}{U-D}g_k(U,D,s_{k-1}U)-g_k(U,D,s_{k-1}\xi _k).
\label{eq_residual}
\end{equation}
In this case where $U$ and $D$ are known, and consequently $D\le\xi
_k\le U,$ it follows from the convexity of the pay-off function $f$
that the residual is non-negative:
\begin{equation}
\delta_k(U,D)\geq 0,\quad k=1,\dots ,n.
\end{equation}
 In this fashion, each stock price
process path $\{s_k\}$ maps to a corresponding sequence of
non-negative residuals $\{\delta_k(U,D)\}$, which are withdrawn
after each portfolio liquidation prior to the construction of the
next time period super-hedge. The accumulated value of the
withdrawn residuals at maturity is given by
\begin{equation}
\Delta_n(U,D)=\delta_1(U,D)(1+r)^{n-1}+\delta_2(U,D)(1+r)^{n-2}+\cdots
+\delta_n(U,D). \label{eq_Delta}
\end{equation}

For the remainder of this paper we will assume a European call
option pay-off function $f,$
\begin{equation}
f(s)=(s-K)_+
\end{equation}
where $K$ is the option strike price.
\section{Model Properties}

Ruschendorf \cite{Ruschendorf2002} and Nagaev \cite{Nagaev2003}
document a number of properties of our financial market model which
will illuminate our algorithmic design approach. In particular, the
incompleteness of the market model is manifested in an open interval
$(\underline{x}_k ,\overline{X}_k)$ $(k=0,\dots ,n-1)$ of
no-arbitrage option prices. The end points of the interval are shown
to be
\begin{align*}
    \overline{X}_k(U,D) &= g_k(U,D,s_k) \\
    \underline{x}_k(U,D) &= (1+r)^{-(n-k)} \left(s_k (1+r)^{n-k} -
    K\right)_+
\end{align*}
where $K$ is the option strike price and $g_k$ is defined in
(\ref{CRR_price}).

At every time instant $t = k$ $(k=0,\dots ,n-1)$, the open interval
$(\underline{x}_k(U,D),\overline{X}_k(U,D))$ is the set of
no-arbitrage option prices. For the option seller, the upper bound
$\overline{X}_k(U,D)$ is the demarcation between risk sharing with
the option buyer (if the option sale price is below
$\overline{X}_k(U,D)$) and the potential for arbitrage profit (if
the option sale price is at or above $\overline{X}_k(U,D)$).

Let us now replace the parameters $U,$ $D$ with a pair of numbers
$u,$ $d$ such that
\begin{equation} D \le d < u\le U\label{ineq}\end{equation}
and let us construct the portfolio with the time $t=k$ set-up cost
\begin{equation}\overline{x}_k(u,d)=g_k(u,d,s_k),\label{setup-u-d}\end{equation}
where $g_k$ is defined in (\ref{CRR_price}). It is straightforward
to show that for any choice of $d$ and $u$ satisfying (\ref{ineq}),
the resulting portfolio set-up cost (\ref{setup-u-d}) falls within
the no-arbitrage option price interval
\[
\underline{x}_k(U,D) \le \overline{x}_k(u,d) \le
\overline{X}_k(U,D), k=0,\dots ,n.
\]
We will refer to the value $\overline{x}_k(u,d)$ as a rational
price of the option.
\section{Choice of $(u,d)$ Pair}
\label{section:choice}
 As a practical matter, we do not know the
actual values of $D$ and $U$. Suppose we choose a $(u,d)$ pair that
is within the $D$ and $U$ values:
\[
D \le d < u \le U
\]
and suppose that for any given stock price process $\{s_k\},$ we
define the hedging portfolio strategy
\begin{equation}
(\gamma_k(u,d),\beta_k(u,d)),~k=0,\dots,n-1 \label{eq_portfolio}
\end{equation}
where $\gamma_k(u,d)$ and $\beta_k(u,d)$ are defined in
(\ref{eq_gamma}) and (\ref{eq_beta}) respectively, with the boundary
parameters $U,D$ replaced with the values $u,d$. The above portfolio
strategy will produce a residual sequence
\begin{equation}
\delta_k(u,d),~k=1,\dots,n \label{eq_resid2}
 \end{equation}
 where
$\delta_k(u,d)$ is defined in (\ref{eq_residual}), with $U,D$
replaced by $u,d$. It is straightforward to show that
\begin{itemize}
\item $\delta _k(u,d)>0$  if $d < \xi _k < u$
\item $\delta _k(u,d)=0$  if $\xi _k=d$ or $\xi _k=u$
\item $\delta _k(u,d)<0$  if $D < \xi _k<d$ or $u < \xi _k < U.$
\end{itemize}
In order to maintain the dynamic portfolio strategy defined by
(\ref{eq_portfolio}), at each time step $k=1,\dots ,n$ the investor
will withdraw the residual (\ref{eq_resid2})  from the liquidated
proceeds when $\delta_k(u,d)>0$
 and add the amount when
 $\delta_k(u,d)<0$. The risk-free growth of the local residuals $\delta_k(u,d)$ produces an accumulated residual defined by
\begin{equation}
\Delta_n(u,d)=\delta_1(u,d)(1+r)^{n-1}+\delta_2(u,d)(1+r)^{n-2}+\cdots
+\delta_n(u,d)\label{acc_resid2}.
\end{equation}

We would like to stress here that the hedging portfolio strategy
constructed above is in general non-self-financing. An investor who
utilizes our dynamic hedging strategy will want to choose
 values for $d$ and $u$ that determine a residual sequence with
desirable statistical characteristics. It is the choice of the
model parameter values $d$ and $u$ based on the statistical
characteristics of the residual sequence that constitutes our
algorithm design.

\section{Risk Minimization on $(u,d)$ Contours}
\label{sec:riskmin}
 We propose a two-stage algorithm for choosing a
$(u,d)$ pair. The first stage reduces the set of $(u,d)$ pairs under
consideration by imposing a market calibration constraint. The
second stage chooses from this reduced set a pair that optimizes one
of a number of investor-relevant statistical properties of the
residual sequence.

\subsection{Market Calibrated Price Contour}
The first stage of our proposed risk minimization procedure is the
selection of a set of $(u,d)$ pairs consistent with the quoted
market option price.

Each $(u,d)$ pair uniquely determines a portfolio strategy that is
dependent upon the realized values of the stock and bond processes.
At the initial time $t=0$, the portfolio strategy determined by a
$(u,d)$ pair specifies an initial portfolio consisting of
$\beta_0(u,d)$ bonds and $\gamma_0(u,d)$ stocks where $\beta_0(u,d)$
and $\gamma_0(u,d)$ are defined in (\ref{eq_beta}) and
(\ref{eq_gamma}) with $k=0,~U=u$ and $D=d.$ At time $t=0$ the set-up
cost of the so constructed portfolio is
\begin{equation}
g_0(u,d,s_0) = (1+r)^{-n}\sum_{j=0}^{n}C_{n}^j[p(u,d)]^j
[1-p(u,d)]^{n-j}(s_0u^jd^{n-j}-K)_+
\end{equation}
where $p(u,d)$ is given by (\ref{eq_pval}) (with $U,D$ replaced by
$u,d$) and $s_0$ is the stock price at $t=0.$ To calibrate the
choice of $(u,d)$ pairs to the market price of the option, we need
to set $g_0(u,d,s_0)$ equal to the option time $t=0$ market price
$x_0:$
$$g_0(u,d,s_0)=x_0.$$

\begin{figure}[!h]
\begin{center}
\includegraphics[scale=.5]{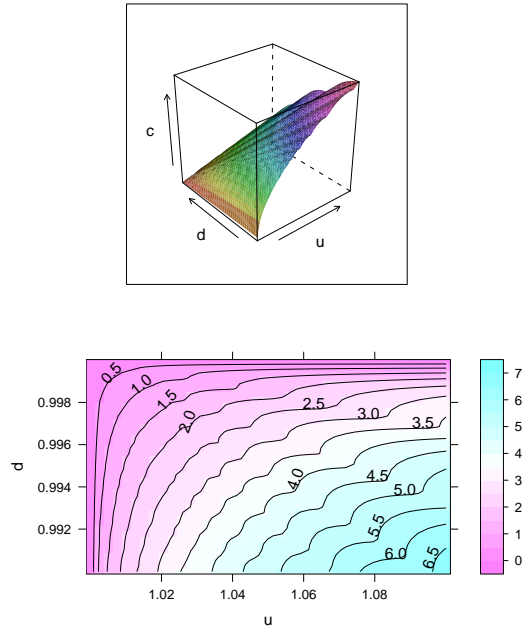}
\end{center}
\caption{Value Surface and Contours} \label{fig:surface}
\end{figure}

While the choice of $(u,d)$ uniquely determines a portfolio set-up
cost, specifying a portfolio set-up cost determines a contour of
$(u,d)$ pairs since there are a multiplicity of portfolios with
identical set-up costs. This contour consists of the set

\begin{equation}
 \Sigma=\{(u,d): c_0(u,d)=c^* \},  \qquad c^* = \frac{x_0}{s_0}
 \label{eq_Sigma}
\end{equation}
where $c_0$ is the normalized value surface
\begin{equation}
  c_0(u,d)=\frac{g_0(u,d,s_0)}{s_0}=\sum_{j=0}^n C^j_n [p(u,d)]^j[1-p(u,d)]^{n-j}(u^jd^{n-j}-R)_+, \qquad R=\frac{K}{s_0}
  \label{eq:surface}
 \end{equation}
(see Figure \ref{fig:surface}).

Computationally, we utilize contour construction software to compute
a finite number of $(u,d)$ pairs satisfying (\ref{eq_Sigma}). It is
this set of $(u,d)$ pairs that is used by the second stage of our
algorithm.

\subsection{Investor-relevant Choice Criteria}
\label{section:criteria} The second stage of the two-stage algorithm
selects a unique $(u,d)$ pair on the market calibrated contour
$\Sigma$ defined in (\ref{eq_Sigma}) that optimizes one of several
investor risk criteria.
 Each $(u,d)$ pair on the market-calibrated
contour $\Sigma$ determines a \cor{dynamic} hedging portfolio
strategy. For a given stock price process $\{s_k\}$, the
$(u,d)$-determined portfolio\cor{s} produce a sequence of residuals
$\delta_k(u,d),$  each representing a residual profit/loss for an
investor (see section \ref{section:choice} for details). This
sequence is an economic measure of the consequence of choosing model
parameters $(u,d)$ and the associated dynamic hedging portfolio.

\cor{There are s}everal criteria that convert this sequence into a
scalar measure of investor risk, \cor{each reflecting some aspect
of the option seller attitude towards risk. We thus have the
following situation. Fix a $(u,d)$ pair at time $t=0$. Each
potential stock price time series $\{s_k\},~k=1,\cdots,n$
determines a sequence of residuals
$\{\delta_k(u,d),~k=1,\cdots,n\}$. A particular choice of a risk
criterion reduces the sequence $\{\delta_k(u,d)\}$ to a single
scalar value of risk. To judge the acceptability of the $(u,d)$
pair under the chosen risk measure, we simulate a number of stock
price time series and collect the corresponding sample of scalar
risk values. An appropriate sample statistic (mean value or
probability of a desirable event) is then computed from the sample
as the utility value of the $(u,d)$ pair. The behavior of the
sample statistic as the $(u,d)$ pair is varied over the
market-calibrated contour determines the optimal choice of
$(u,d)$.}

There are four criteria that we consider for choosing a unique
$(u,d)$ pair.

\begin{itemize}
\item Maximize the likelihood of a positive accumulated residual:
\begin{equation}
\max_{(u,d)\in ~\Sigma} \mbox{prob}(\Delta_n(u,d) > 0) \label{opt1}
\end{equation}

\item Minimize expected shortfall:
\begin{equation}
\min_{(u,d)\in~\Sigma}E(shortfall) =
\min_{(u,d)\in~\Sigma}E(\max_{1 \le k \le n}(-\delta_k(u,d)))
\label{opt2}
\end{equation}

\item Minimize the expected accumulated squared residuals:
\begin{equation}
\min_{(u,d)\in~\Sigma} E(\sum_{k=1}^n (\delta_k(u,d))^2)
\label{opt3}
\end{equation}

\item Maximize the expected accumulated profit
\begin{equation}
\max_{(u,d)\in~\Sigma} E(\Delta_n(u,d)) \label{opt4}
\end{equation}
\end{itemize}

The first criteri\cor{on}  (\ref{opt1}) interprets a positive
residual as a profit, and chooses a $(u,d)$ pair that has the
highest probability of a net profit. In the absence of arbitrage a
large accumulated profit is not attainable with high probability.
There is the possibility, however, of an investor achieving a
small positive profit. The optimization problem presented here
produces a market calibrated hedging strategy that maximizes the
likelihood of a positive accumulated profit.

The second \cor{criterion}  (\ref{opt2}) reflects an investor's
desire to minimize the amount of single period additional funding
needed to rebalance the portfolio over the life of the option.
\cor{A negative}  residual $\delta_k(u,d)$ represents the cash
shortfall of the portfolio value at time $k.$ The largest negative
$\delta_k(u,d)$ is the largest shortfall value. \cor{Optimizing this
criterion} produces a hedging portfolio with minimal expected single
period additional funding.

In (\ref{opt3}) the residual $\delta_k(u,d)$ represents an
economic measure of  model error.  Concern for minimizing model
risk would motivate weighting equally positive and negative
$\delta_k(u,d)$ residuals. \cor{Indifference to the sign of
$\delta_k(u,d)$ can be achieved by using squared residuals in the
risk criterion.}

 Our final criteri\cor{on}
(\ref{opt4}) maximizes the expected accumulated residual, which
reflects total net profit from using the dynamic portfolio based
on the chosen $(u,d)$. It was shown in \cite{Nagaev2003} that the
expected accumulated profit is asymptotically constant on contours
of constant  rational price. We thus anticipate minimal
differences in the expected accumulated profit at each $(u,d)$
pair on our constant rational price contours when $n$ is large.
For small $n,$ empirical results show it is possible to \cor{have
a market contour with non-constant expected accumulated profit}.


In summary, our model provides a measure of the economic impact of
market incompleteness \cor{by the construction of the} residuals
$\delta_k(u,d)$. The criteria presented here \cor{can be
optimized} to determine hedging strategies \cor{that mitigate}
investor risk.

\section{Algorithm Implementation and Numerical Results}
\label{sec:algandresults} In this section, we present an
implementation of our two-stage algorithm for producing an optimal
hedging strategy as measured by the risk criteria presented in
section \ref{section:criteria}. Following this is a presentation of
illustrative results from applying the algorithm to real market
data.

\subsection{Two-stage Algorithm}
\label{section:algo} The choice of an optimal $(u,d)$ pair
proceeds in two stages.
\begin{itemize}
\item Reduction by market calibration of the population of $(u,d)$
pairs to a contour given by (\ref{eq_Sigma}) corresponding to the
quoted market option price.
\item Imposition of a ranking criterion, based on a bootstrap estimated
statistic of the residual sequence.
\end{itemize}

 The numerical procedures to be
described are implemented in the R computer language
(\cite{Team2005}). The flow of computation proceeds as follows.

\begin{enumerate}
\item The contour creation function in R is applied to the normalized value surface given in
(\ref{eq:surface}).
 Contours defined by (\ref{eq_Sigma}) with $r=0$ are identified as depicted in Figure \ref{fig:surface}.
 The contour matching the market option price $x_0$ is chosen. The
software typically identifies between 80 and 100 $(u,d)$ pairs on
the market calibrated contour.
\item Historical daily stock price data is used to create a sequence
of daily stock price jumps. The jumps are separated into groups by
day count of successive price data: next day jumps (e.g. Monday to
Tuesday price jump) and weekend jumps (i.e. Friday to Monday price
jump) constitute the majority of jumps. There are a few single day
mid-week holiday jumps and long weekend jumps. A typical stock
price process and jump process are depicted in Figure
\ref{fig:intel}.
\begin{figure}[!h
]
\begin{center}
\subfigure[Stock Price Time Series]
{\includegraphics[scale=.5]{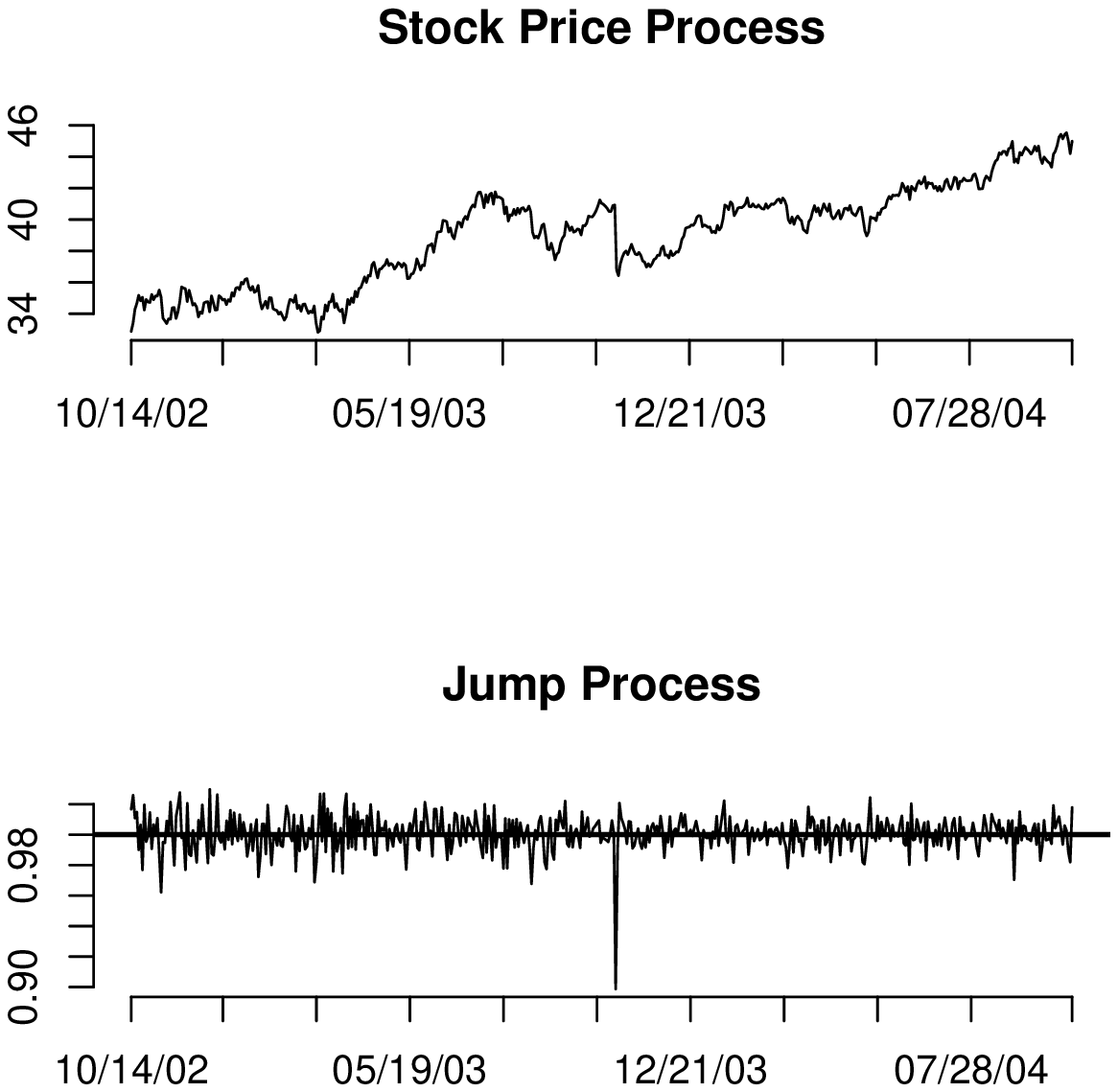} } \qquad
\subfigure[Jump Distributions]
{\includegraphics[scale=.5]{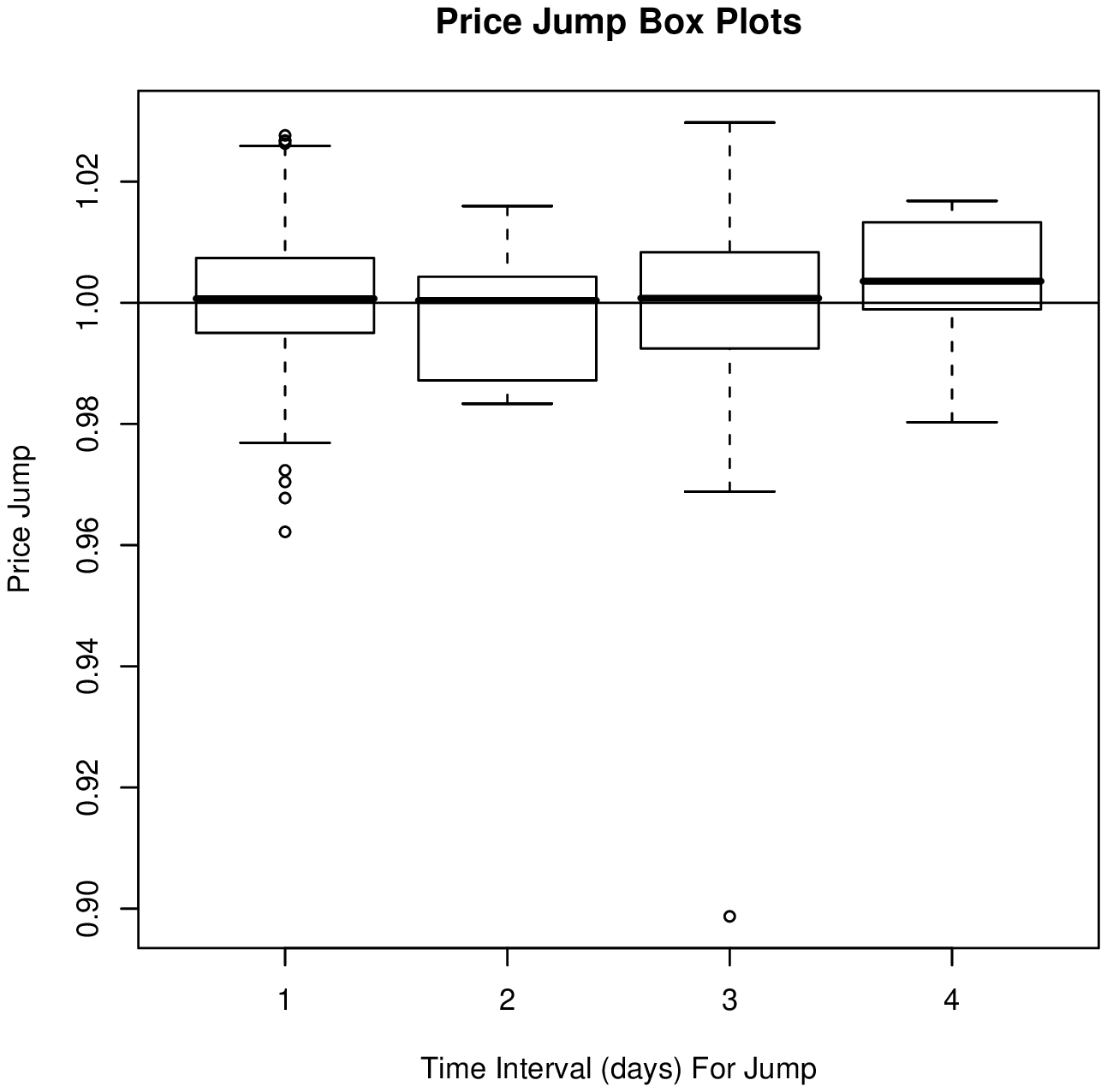} }
\end{center}
\caption{Stock Price Process and Price Jump Process}
\label{fig:intel}
\end{figure}
\item Each of the next day and weekend groups are sampled with replacement
to form bootstrap jump sequences (four next-day jumps followed by
one weekend jump, repeated for however many weeks of bootstrap
data are needed). Each jump sequence uniquely determines a stock
price sequence. A set of bootstrap jump and price processes are
depicted in Figure \ref{fig:bootStocks}.

\begin{figure}[!h]
\begin{center}
\includegraphics[scale=.4]{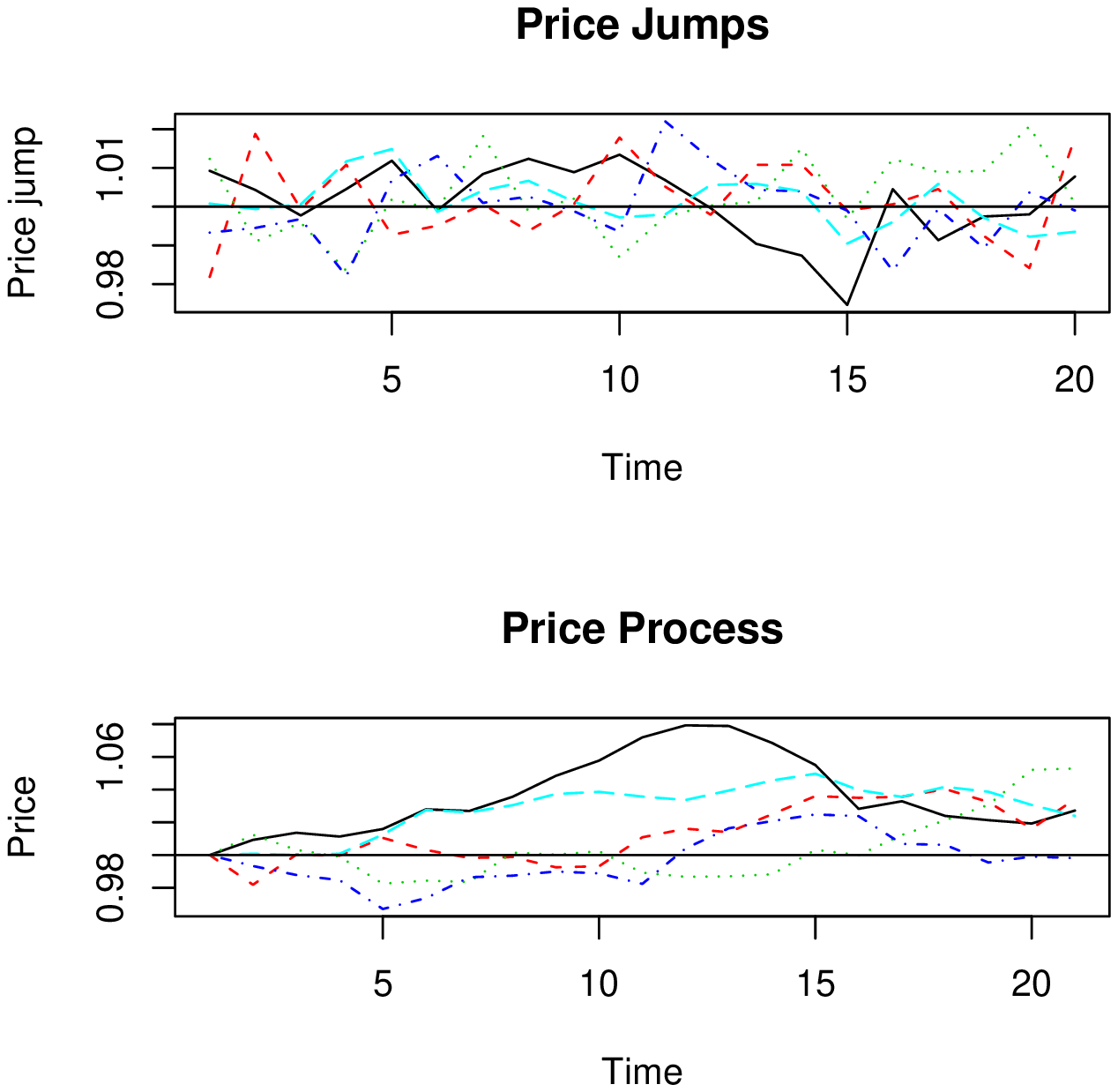}
\end{center}
\caption{Bootstrap Stock Price Processes} \label{fig:bootStocks}
\end{figure}

\begin{figure}[!hb]
\begin{center}
\includegraphics[scale=.5]{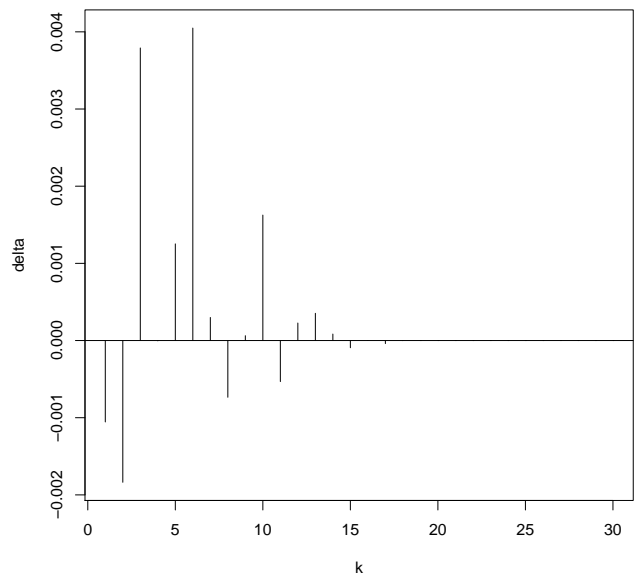}
\end{center}
\caption{Local residual ($\delta_k$) distribution}
\label{fig:resboot}
\end{figure}

\begin{figure}[!hb]
\begin{center}
\includegraphics[scale=.5]{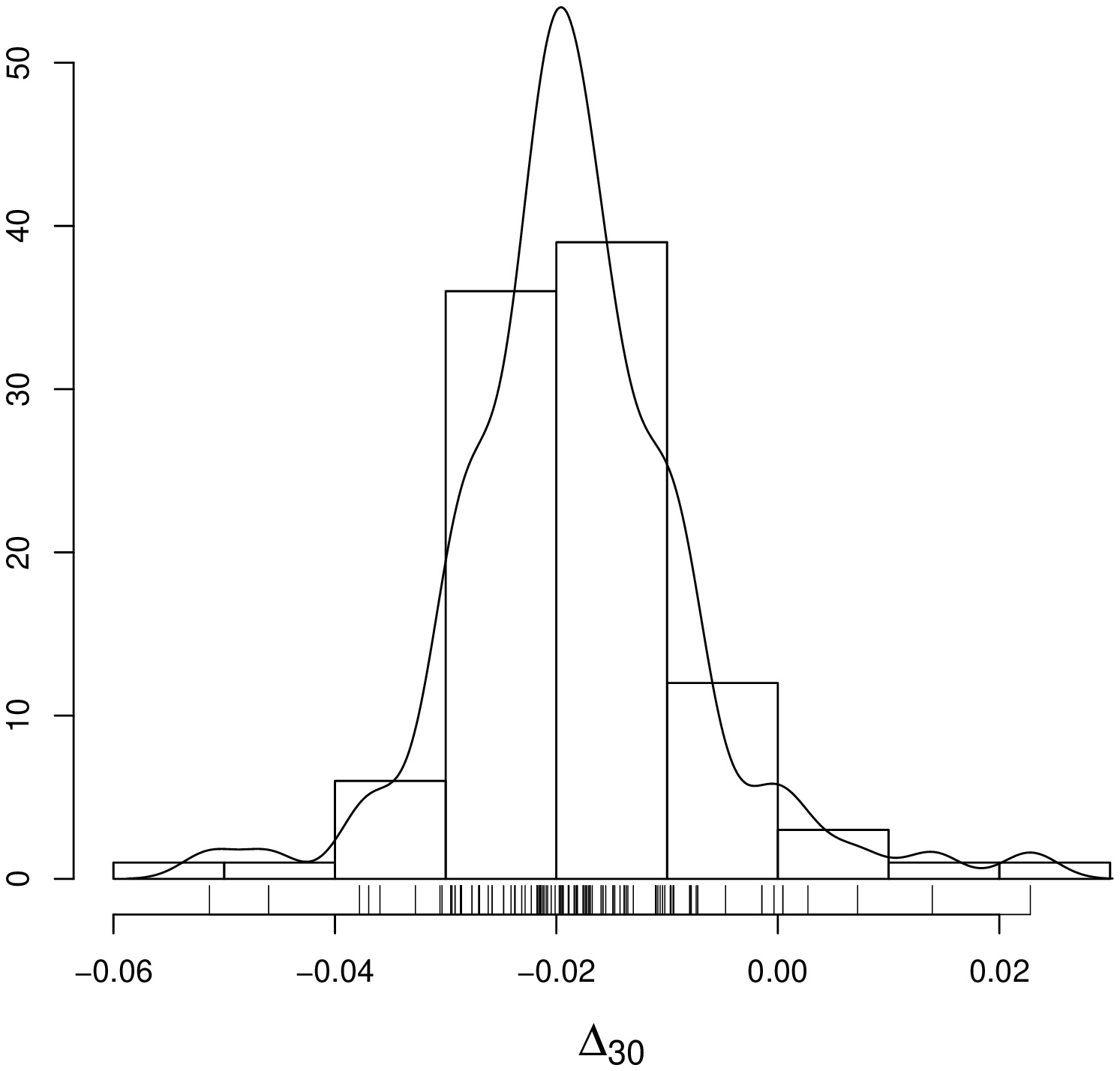}
\end{center}
\caption{Aggregated Residual ($\Delta_n$) Distribution}
\label{fig:bootResiduals}
\end{figure}

\item For each $(u,d)$ pair on the chosen contour, the sequence of
residuals $\delta_k$ are computed for each bootstrap stock price
sequence. A residual sequence for a bootstrap sample with 30 days
to expiration is shown in Figure \ref{fig:resboot}. The
distribution of the corresponding aggregated residual is given in
Figure \ref{fig:bootResiduals}.
\item Each criteria is applied to the bootstrap sample of residual sequences and
the appropriate statistic (expected value or probability) is
estimated.
\item The $(u,d)$ pair with the best criterion value is chosen. This
results in four optimal $(u,d)$ pairs, one for each of the four
criteria described previously.

\end{enumerate}

\subsection{Illustrative Results}
In this section we present results representative of the insight
obtained from applying the risk evaluation tools developed in
previous sections. Results are shown for several call options
traded on the New York Stock Exchange expiring on October 15,
2004. The selected options are as follows
\begin{itemize}
\item Walmart corporation option WMT+JJ, strike price \$50
\item Exxon Mobil corporation option XOM+JV, strike price \$42.50
\item Intel corporation option INQ+JE strike, price \$25.
\end{itemize}
Data collected included historical stock price data, option strike
price and market option price. Over 500 daily stock prices were
recorded and used in constructing bootstrap samples of the stock
price jump process. The market option price data for a period of
40 days prior to expiration was collected and used in selecting
market calibrated contours (as described in section
\ref{section:algo}) with varying time to expiration. The market
option price $n$ days prior to expiration was used in identifying
the appropriate contour for each of the reported values of $n.$

Selecting a market calibrated contour produces a set of
approximately 80 $(u,d)$ pairs each corresponding to a potential
hedging portfolio. The risk criteria (\ref{opt1}) through
(\ref{opt4}) are evaluated for each hedging portfolio and the
optimal is chosen. Table \ref{xom} presents results for the
Exxon-Mobil option with $n=30$ days to expiration. The algorithm
produces a hedging portfolio where the probability of a positive
profit is 1. This is consistent with the theoretical results
presented in section \ref{section:choice}. When the sequence of
relative stock price jumps $\xi_k$ fall within the interval
$(u,d),$ we are guaranteed the local residual profits $\delta_k$
will be positive. The optimal expected shortfall and expected
squared deviation are both essentially zero. For this particular
data the objective values along the contour did not vary
dramatically. The advantage of utilizing the algorithm is seen in
the expected shortfall computation where the optimal expected
shortfall is approximately 3 times as small as other possible
shortfall values.

\begin{table}
\centering
\begin{tabular}{|c|c|c|}
\hline Risk criterion &Optimal (u,d)& Objective value\\
\hline
$P(\Delta_n)>0$&(1.0112,0.9792) &  1.0 \\
\hline
 $E(\mbox{shortfall})$ &(1.0146,0.9834) &  -0.0042\\
     \hline
 $E(\sum_{k=1}^n\delta_k^2)$&(1.0130,0.9822)& 0.0053\\
     \hline
 $E(\Delta_n)$&(1.0100,0.9961) & 0.2066\\
    \hline
\end{tabular}
\caption{Exxon-Mobil optimal results with $n=30$ days to
expiration}
\label{xom}
\end{table}

Results for the Walmart option are presented in tables \ref{wmt30}
through \ref{wmt7}. Examining the results we see that for $n=30$
and $n=20$ days to expiration, the optimal $(u,d)$ pair produces
with probability $1.0$ a small positive aggregated profit. It is
interesting to note that over all $(u,d)$ pairs on the contour,
the probability of a positive aggregated profit ranged from
approximately 0.6 to 1.0 when $n=30$ and $n=20$ and from
approximately 0.3 to 0.99 when $n=7.$ In other words, an investor
can increase the probability of achieving a small positive profit
from 0.3 to 0.99 by following the hedging strategy produced by the
algorithm.

Similar results were seen in comparing the values of the other
risk criteria for the optimal hedging portfolio as compared to
other hedging portfolios associated with $(u,d)$ pairs on the
contour. In particular, considering the Walmart data with $n=30$
days to expiration, the values for the expected shortfall $E(\max
(-\delta_k(u,d)))$ ranged from -0.0548 (see the second column of
\ref{wmt30}) to approximately -0.4. The shortfall is more than 7
times as large as the optimal for some hedging portfolios. The
values of the expected squared deviation $E((\delta_k(u,d))^2)$
ranged from the minimum of 0.0225 to a maximum of approximately
0.18 and the expected aggregated profit $E(\Delta_n(u,d))$ ranged
from approximately 0 to 0.3. These results illustrate the value of
following the hedging strategy suggested by the algorithm.

\begin{table}
\centering
\begin{tabular}{|c|c|c|}
\hline Risk criterion &Optimal (u,d)& Objective value\\
\hline
$P(\Delta_n)>0$&(1.0156,0.9897) &  1.0 \\
\hline
 $E(\mbox{shortfall})$ &(1.0140,0.9889) &  -0.0548\\
     \hline
 $E(\sum_{k=1}^n\delta_k^2)$&(1.0134,0.9881)& 0.0225\\
     \hline
 $E(\Delta_n)$&(1.0076,0.9804) & 0.3097\\
    \hline
\end{tabular}
\caption{Walmart optimal results $n=30$ days to expiration}
\label{wmt30}
\end{table}

\begin{table}
\centering
\begin{tabular}{|c|c|c|}
\hline Risk criterion &Optimal (u,d)& Objective value\\
\hline
$P(\Delta_n)>0$&(1.0171,0.9901) &  1.0 \\
\hline
 $E(\mbox{shortfall})$ &(1.0136,0.9878) &  -0.0355\\
     \hline
 $E(\sum_{k=1}^n\delta_k^2)$&(1.0136,0.9878)& 0.0220\\
     \hline
 $E(\Delta_n)$&(1.0052,0.9718) & 0.2372\\
    \hline
\end{tabular}
\caption{Walmart optimal results $n=20$ days to expiration}
\label{wmt20}
\end{table}

\begin{table}
\centering
\begin{tabular}{|c|c|c|}
\hline Risk criterion &Optimal (u,d)& Objective value\\
\hline
$P(\Delta_n)>0$&(1.0238,0.9882) &  0.99 \\
\hline
 $E(\mbox{shortfall})$ &(1.0178,0.9857) &  -0.0018\\
     \hline
 $E(\sum_{k=1}^n\delta_k^2)$&(1.0116,0.9820)& 0.0054\\
     \hline
 $E(\Delta_n)$&(1.0068,0.9747) & 0.0215\\
    \hline
\end{tabular}
\caption{Walmart optimal results $n=7$ days to expiration}
\label{wmt7}
\end{table}

Comparing tables \ref{wmt30} through \ref{wmt7} we see that the
optimal objective function value for most of the risk criteria
does not change much as the length of time to expiration changes
with the exception of the expected aggregated profit. The value
for $n=30$ is more than 10 times larger than the value for $n=7$.
The difference could be explained by the fact that given more time
to expiration, there are more opportunities to withdraw a small
positive profit (at each $k=1,2,\cdots,30$). If we consider the
results produced for the Intel option for $n=30,~n=20$ and $n=10$
presented in tables \ref{intel30} through \ref{intel10},  we do
not see the same change in the value of the expected aggregated
profit. In this case, however, the probability of achieving a
positive profit is more variable with the length of time to
expiration increasing from 0.8 to 0.93 as $n$ decreases from
$n=30$ to $n=10.$ Clearly, the behavior of the optimal objective
values vary with each chosen option.

\begin{table}
\centering
\begin{tabular}{|c|c|c|}
\hline Risk criterion &Optimal (u,d)& Objective value\\
\hline
$P(\Delta_n)>0$&(1.0140,0.9640) &  0.8 \\
\hline
 $E(\mbox{shortfall})$ &(1.0130,0.9616) &  -0.0378\\
     \hline
 $E(\sum_{k=1}^n\delta_k^2)$&(1.0160,0.9696)& 0.0104\\
     \hline
 $E(\Delta_n)$&(1.0127,0.9608) & 0.0429\\
    \hline
\end{tabular}
\caption{Intel optimal results with $n=30$ days to expiration}
\label{intel30}
\end{table}

\begin{table}
\centering
\begin{tabular}{|c|c|c|}
\hline Risk criterion &Optimal (u,d)& Objective value\\
\hline
$P(\Delta_n)>0$&(1.0193,0.9624) &  0.83 \\
\hline
 $E(\mbox{shortfall})$ &(1.0182,0.9600) &  -0.0098\\
     \hline
 $E(\sum_{k=1}^n\delta_k^2)$&(1.0229,0.9703)& 0.0045\\
     \hline
 $E(\Delta_n)$&(1.0197,0.9632) & 0.0253\\
    \hline
\end{tabular}
\caption{Intel optimal results with $n=20$ days to expiration}
\label{intel20}
\end{table}

\begin{table}
\centering
\begin{tabular}{|c|c|c|}
\hline Risk criterion &Optimal (u,d)& Objective value\\
\hline
$P(\Delta_n)>0$&(1.0339,0.9600) &  0.93 \\
\hline
 $E(\mbox{shortfall})$ &(1.0130,0.9616) &  -0.0034\\
     \hline
 $E(\sum_{k=1}^n\delta_k^2)$&(1.0160,0.9696)& 0.0062\\
     \hline
 $E(\Delta_n)$&(1.0127,0.9608) & 0.0475\\
    \hline
\end{tabular}
\caption{Intel optimal results with $n=10$ days to expiration}
\label{intel10}
\end{table}

\section{Conclusions}
\label{sec:conclusion} We have developed an algorithm that produces
a non-self-financing hedging strategy in an incomplete market
corresponding to one of several investor risk criteria. The
algorithm provides the opportunity to evaluate the economic
consequences of choosing a particular hedging strategy in an
incomplete market. The two-stage algorithm optimizes one of a number
of investor-relevant statistical properties of a local residual
profit or shortfall.

The algorithm was tested on several options traded on the New York
stock exchange. The results illustrate that following the
portfolio strategy produced by the algorithm is beneficial to an
investor, improving the value of the investor risk criterion by as
much as a factor of ten compared to the results associated with
other, non-optimal hedging portfolio strategies.

\section{Future Research}
\label{sec:future}
 In this paper we investigate non-self financing hedging strategies for a short-term
European call option (time to expiration, $n,$ is at most 30
days). Our algorithm builds a short-term portfolio strategy that
optimizes one of the investor-related criteria. Our research was
inspired by theoretical investigations of A. Nagaev et al. (see
\cite{Nagaev2003}, \cite{WiennaMarkov},
\cite{WiennaNon-Identical}, \cite{SashaRussian}) where the
long-term behavior (large $n$) of the accumulated residual
(\ref{acc_resid2})
 has been studied.

In \cite{Nagaev2003}, \cite{WiennaMarkov}, and
\cite{WiennaNon-Identical}, asymptotic properties of the so-called
riskless profit of an investor ($\Delta _n(U,D)$ defined in
(\ref{eq_Delta}) with the boundary parameters $U$ and $D$) have
been studied. The  case of  independent identically distributed
(i.i.d.) stock price jumps is presented in \cite{Nagaev2003}; in
\cite{WiennaMarkov} the jumps are assumed to follow a discrete
Markov chain; \cite{WiennaNon-Identical} studies the case of
independent, but not identically distributed jumps. In all three
cases, by means of suitable diffusion approximations, asymptotic
formulas for the mean accumulated residuals in terms of the
original model parameters have been obtained.

In \cite{SashaRussian}, A. Nagaev considers the accumulated residual
(\ref{acc_resid2}) for the case of non-boundary parameters $u$ and
$d$ satisfying (\ref{ineq}) (the so-called risky profit of an
investor). In this case, under the i.i.d. assumptions on the stock
price jumps, the asymptotic formula for the mean accumulated
residual has been obtained and asymptotic connections (as time to
expiration $n$ tends to infinity) between $E(\Delta _n)$ and the
set-up cost $g_0(u,d,s_0)$ have been established.

The present paper is an initial study of the short-term
accumulated residuals. We assume here that the relative stock
price jumps are i.i.d. random variables and build our bootstrap
simulation procedure accordingly. In the future we will
investigate the consequences of more realistic assumptions on
stock price jumps using methods of time series analysis and/or
advanced model fitting (e.g. Levy processes based probability
models). We also plan to extend our investigations to other
derivative securities with convex pay-off functions. An additional
direction of future research is numerical testing of the
asymptotic formulas obtained by A. Nagaev et al., more
specifically, exploring the consequences of the finite number of
time steps to expiration on asymptotic results.

\section{Acknowledgements}
We would like to dedicate this paper to our colleague, Alexander
Nagaev. Alexander was an insightful mathematician, creative thinker
and a true friend. He died unexpectedly in 2005 and we truly miss
him.

This work is partially supported by Bentley College Enterprise Risk
Management Research Grants.

\bibliographystyle{ieeetran}
\bibliography{lucy}

\end{document}